\documentclass[14pt]{amsart}
\address{\newline{\normalsize University of Edinburgh, Kings Buildings, Mayfield
Road, Edinburgh EH9 3JZ, UK}\newline{\it E-mail address}:
ilkarjem@rambler.ru}
\usepackage{amscd,amssymb}
\usepackage{amsthm,amsmath,amssymb}
\usepackage[matrix,arrow]{xy}

\newtheorem{theorem}[equation]{Theorem}
\newtheorem{proposition}[equation]{Proposition}
\newtheorem{lemma}[equation]{Lemma}
\newtheorem{corollary}[equation]{Corollary}
\newtheorem{conjecture}[equation]{Conjecture}

\theoremstyle{definition}

\theoremstyle{remark}
\newtheorem{remark}[equation]{Remark}

\makeatletter\@addtoreset{equation}{section}\makeatother

\begin{document}


\title{Remark on polarized K3 surfaces of genus 36}

\thanks{The work was partially supported by
RFFI grant No. 08-01-00395-a and grant N.Sh.-1987.2008.1.}

\author{Ilya Karzhemanov}

\begin{abstract}
Smooth primitively polarized $\mathrm{K3}$ surfaces of genus $36$
are studied. It is proved that all such surfaces $S$, for which
there exists an embedding $\mathrm{R} \hookrightarrow
\mathrm{Pic}(S)$ of some special lattice $\mathrm{R}$ of rank $2$,
are parameterized up to an isomorphism by some $18$-dimensional
unirational algebraic variety. More precisely, it is shown that a
general $S$ is an anticanonical section of a (unique) Fano
$3$-fold with canonical Gorenstein singularities.
\end{abstract}

\sloppy

\maketitle

\section{Introduction}
\label{section:introduction}

Let $\mathcal{K}_{g}$ be the moduli space of all smooth
primitively polarized $\mathrm{K3}$ surfaces of genus $g$.
$\mathcal{K}_{g}$ is known to be a quasi-projective algebraic
variety (see for example \cite{viehweg}). This makes it possible
to consider the fundamental questions of birational geometry about
$\mathcal{K}_{g}$ such as its rationality, unirationality,
rational connectedness, Kodaira dimension, and etc.

S. Mukai's vector bundle method, developed to classify higher
dimensional Fano manifolds of Picard number $1$ and coindex $3$
(see \cite{mukai-1}, \cite{mukai-2}), allowed to prove
unirationality of $\mathcal{K}_{g}$ for $g \in \{2, \ldots, 10,
12, 13, 18, 20\}$ (see \cite{mukai-3}, \cite{mukai-4},
\cite{mukai-5}, \cite{mukai-6}). At the same time,
$\mathcal{K}_{g}$ turns out to be non-unirational for general $g
\geqslant 43$ (see \cite{gritsenko-hulek-sancaran},
\cite{kondo-1}, \cite{kondo-2}). In principle, the proof of
unirationality of $\mathcal{K}_{g}$ is based on the observation
that general $\mathrm{K3}$ surface $S_g$ with primitive
polarization $L_g$ and ``not very big'' $g$ is an anticanonical
section of a smooth Fano $3$-fold $X_g$ of genus $g$ so that $L_g
= -K_{X_g}\big\vert_{S_g}$ (see \cite{mukai-3}, \cite{mukai-5},
\cite{mukai}). The latter gives a rational dominant map from the
moduli space $\mathcal{F}_g$ of pairs $(X_{g}, S_{g})$, where
$S_{g} \in |-K_{X_{g}}|$ is smooth, to $\mathcal{K}_{g}$ by
sending $(X_{g}, S_{g})$ to $S_{g}$, with $\mathcal{F}_g$
typically being a rational algebraic variety. However, this
construction has the restriction that $X_g$ must have Picard
number $1$, which does not hold for most $g$ (see
\cite{isk-prok}).

In order to generalize the above arguments for every possible $g$,
to a given smooth Fano $3$-fold $V$ of genus $g$ one associates
the Picard lattice $R_V := \mathrm{Pic}(V)$, equipped with the
pairing $(D_{1}, D_{2}) := D_{1} \cdot D_{2} \cdot (-K_{V})$ for
$D_1$, $D_2 \in \mathrm{Pic}(V)$, and considers the moduli space
$\mathcal{K}_{g}^{R_{V}}$ of all smooth $\mathrm{K3}$ surfaces
$S_g$, equipped with a primitive embedding $R_V \hookrightarrow
\mathrm{Pic}(S_{g})$ which maps $-K_V$ to an ample class on $S_g$
of square $g$ (let us call such $S_g$ a \emph{$\mathrm{K3}$
surface of type $R_{V}$}). A beautiful result due to A. Beauville
states that a general $\mathrm{K3}$ surface of type $R_{V}$ is the
anticanonical section of a smooth Fano $3$-fold $X_g$ of genus $g$
such that $R_{X_g} \simeq R_{V}$ (see \cite{beauville}). The proof
employs the same idea as above, but instead of $\mathcal{F}_{g}$
the moduli space $\mathcal{F}_{g}^{R_V}$ of pairs $(X_{g},
S_{g})$, where $S_{g} \in |-K_{X_{g}}|$ is smooth and $X_g$ is
equipped with the lattice isomorphism $R_{X_g} \simeq R_V$, is
considered. Again the forgetful map $(X_{g}, S_{g}) \mapsto S_g$
from $\mathcal{F}_{g}^{R_V}$ to $\mathcal{K}_{g}^{R_V}$ turns out
to be generically surjective. However, these arguments can be
applied only to some $g \leqslant 33$ (see \cite{isk-prok}).

In the present paper, we study primitively polarized smooth
$\mathrm{K3}$ surfaces of genus $36$ and consider the following

\begin{conjecture}
\label{theorem:conjecture} The moduli space $\mathcal{K}_{36}$ is
unirational.
\end{conjecture}

To develop an approach to prove
Conjecture~\ref{theorem:conjecture}, we employ the above ideas to
realize a general smooth primitively polarized $\mathrm{K3}$
surface of genus $36$ as an anticanonical section of some Fano
$3$-fold, which must be singular in this case (see
\cite{isk-prok}). The natural candidate for the latter is the Fano
$3$-fold $X$ with canonical Gorenstein singularities and genus
$36$, constructed and studied in \cite{karz-1}, \cite{karz-2}.
This $X$ has only one singular point (see
Corollary~\ref{theorem:singularities-of-x}) and the anticanonical
linear system $|-K_{X}|$ gives an embedding $X \hookrightarrow
\mathbb{P}^{37}$ (see
Remark~\ref{remark:x-is-anti-canonically-embedded}), which implies
that a general surface $S \in |-K_{X}|$ is smooth. Also the Picard
group of $X$ is generated by $K_{X}$ (see
Corollary~\ref{theorem:k-x-generates-pic-of-x}).

Unfortunately, the divisor class group of $X$ has two generators,
$K_X$ and some surface $E$ (see
Corollary~\ref{theorem:k-x-generates-pic-of-x}), so that the
restrictions $K_{X} \big\vert_{S}$ and $E \big\vert_{S}$ generate
a primite sublattice $R_S$ in $\mathrm{Pic}(S)$. In particular,
the Picard number of $S$ must be at least $2$, and hence $S$ can
not be general. However, all lattices $R_S$, $S \in |-K_{X}|$, are
isomorphic to the lattice $\mathrm{R} \simeq \mathbb{Z}^2$ with
the associated quadratic form $70x^2 + 4xy - 2y^2$ (see the end of
Section~\ref{section:preliminaries}), and, as above, we can
consider the moduli space $\mathcal{K}_{36}^{\mathrm{R}}$ of
$\mathrm{K3}$ surfaces of type $\mathrm{R}$. On the other hand, we
may also consider the moduli space $\mathcal{F}$ of pairs
$(X^{\sharp}, S^{\sharp})$, where $X^{\sharp}$ is a Fano $3$-fold
of genus $36$ with canonical Gorenstein singularities and
$S^{\sharp} \in |-K_{X^{\sharp}}|$ is smooth (see
Remark~\ref{remark:proof-of-main-theorem} below for the precise
description of $\mathcal{F}$). Let us state the main result of the
present paper:

\begin{theorem}
\label{theorem:main} The forgetful map $s : \mathcal{F}
\longrightarrow \mathcal{K}_{36}^{\mathrm{R}}$ is generically
surjective.
\end{theorem}

\begin{remark}
\label{remark:proof-of-main-theorem} In the proof of
Theorem~\ref{theorem:main}, we do not appeal to Akizuki--Nakano
Vanishing Theorem, used in \cite{beauville} to show that
$\mathcal{F}_{g}$ (or $\mathcal{F}_{g}^{R_V}$) is a smooth stack,
since it is not clear how to apply this theorem in the singular
case. Instead, we note that $X$ is unique up to an isomorphism
(see Proposition~\ref{theorem:unique-up-to-isomorphism}), and,
moreover, it admits a crepant resolution $f : Y \longrightarrow
X$, with $Y$ being also unique up to an isomorphism (see
Proposition~\ref{theorem:uniqueness-up-to-isomorphism}). Then one
can prove (see Proposition~\ref{theorem:locally-free-action}) that
$\mathcal{F}$ carries the structure of a normal scheme, being the
geometric quotient $U/\mathrm{Aut}(Y)$ of an open subset $U$ in
$\mathbb{P}^{37}$ by the group $\mathrm{Aut}(Y)$ of regular
automorphisms of $Y$. The proof of Theorem~\ref{theorem:main} then
goes along the same lines as in \cite{beauville} (see
Lemma~\ref{theorem:s-is-surjective} below).
\end{remark}

\begin{remark}
\label{remark:proof-of-main-theorem-1} Taking $X =
\mathbb{P}(1,1,1,3)$ in the above considerations, one might apply
the arguments from \cite{beauville} directly (cf.
Remark~\ref{remark:proof-of-main-theorem}) to prove that the
moduli space $\mathcal{K}_{10}$ is unirational (see \cite{karz-1},
\cite{karz-2} for geometric properties of $\mathbb{P}(1,1,1,3)$).
\end{remark}

Furthermore, since the forgetful map
$\mathcal{K}_{36}^{\mathrm{R}} \longrightarrow \mathcal{K}_{36}$
is finite and representable (see \cite[(2.5)]{beauville}), from
Theorem~\ref{theorem:main}, construction of $\mathcal{F}$ and
quasi-projectivity of $\mathcal{K}_{36}$ we deduce the following

\begin{corollary}
\label{theorem:unirationality} There exists a $18$-dimensional
unirational algebraic variety which parameterizes up to an
isomorphism all smooth $\mathrm{K3}$ surfaces of type
$\mathrm{R}$. For general such surface $S$, $S \in |-K_{X}|$ and
the Picard lattice of $S$ is isomorphic to $\mathrm{R}$.
\end{corollary}

\begin{remark}
\label{remark:anticanonical-section} On the opposite, it follows
from the proof of Theorem~\ref{theorem:main} and
\cite{cheltsov-1}, \cite{cheltsov-2}, \cite{prokhorov-enriques}
that no general smooth primitively polarized $\mathrm{K3}$ surface
$S$ of genus $36$ can be an ample anticanonical section of a
normal algebraic $3$-fold, except for the cone over $S$.
\end{remark}

\begin{remark}
\label{remark:1-to-1-rationality-etc}
Corollary~\ref{theorem:unirationality} gives only unirational
hypersurface in $\mathcal{K}_{36}$ but not the whole
$\mathcal{K}_{36}$, and hence the proof of
Conjecture~\ref{theorem:conjecture} is still to go. It would be
also interesting to know whether the map $s$ from
Theorem~\ref{theorem:main} is $1$-to-$1$ and
$\mathcal{K}_{36}^{\mathrm{R}}$ is rational (it follows from the
proof of Theorem~\ref{theorem:main} that $s$ is generically
\'etale).
\end{remark}

\bigskip

I would like to thank Yu. G. Prokhorov for drawing my attention to
the problem and many helpful discussions. I also would like to
thank G. Brown, I. Cheltsov, S. Galkin, A. Veselov for fruitful
conversations.

\section{Notation and conventions}
\label{section:notation-and-conventions}

We use standard notions and facts from the theory of minimal
models (see \cite{kollar-mori}, \cite{kawamata-matsuda-matsuki}).
We also use standard notions and facts from the theory of
algebraic varieties and schemes (see \cite{hartshorne-ag}). All
algebraic varieties are assumed to be defined over $\mathbb{C}$.
Throughout the paper we use standard notions and notation from
\cite{kollar-mori}, \cite{kawamata-matsuda-matsuki},
\cite{hartshorne-ag}. However, let us introduce some:

\begin{itemize}

\item We denote by $\mathrm{Sing}(V)$ the singular locus of an algebraic variety
$V$. For $P \in \mathrm{Sing}(V)$, we denote by $(O \in V)$ the
analytic germ of $V$ at $P$.

\smallskip

\item For a $\mathbb{Q}$-Cartier divisor $M$ and an algebraic cycle $Z$ on a normal
algebraic variety $V$, we denote by $M \big \vert_{Z}$ the
restriction of $M$ to $Z$. We denote by $Z_{1} \cdot \ldots \cdot
Z_{k}$ the intersection of algebraic cycles $Z_{1}, \ldots,
Z_{k}$, $k \in \mathbb{N}$, in the Chow ring of $V$.

\smallskip

\item $M_{1}
\equiv M_{2}$ (respectively, $Z_{1} \equiv Z_{2}$) stands for the
numerical equivalence of two $\mathbb{Q}$-Cartier divisors
$M_{1}$, $M_{2}$ (respectively, of two algebraic $1$-cycles
$Z_{1}$, $Z_{2}$) on a normal algebraic variety $V$. We denote by
$\rho(V)$ the Picard number of $V$. $D_1 \sim D_2$ stands for the
the linear equivalence of two Weil divisors $D_1$, $D_2$ on $V$.
We denote by $N_{1}(V)$ the group of classes of algebraic cycles
on $V$ modulo numerical equivalence. We denote by $\mathrm{Cl}(V)$
(respectively, $\mathrm{Pic}(V)$) the group of Weil (respectively,
Cartier) divisors on $V$ modulo linear equivalence.

\smallskip

\item A normal algebraic three-dimensional variety $V$ is called \emph{Fano threefold}
if it has at worst canonical Gorenstein singularities and the
anticanonical divisor $-K_{V}$ is ample. A normal algebraic
three-dimensional variety $V$ is called \emph{weak Fano threefold}
if it has at worst canonical singularities and the anticanonical
divisor $-K_{V}$ is nef and big. The number $(-K_{V})^{3}$
(respectively, $\displaystyle\frac{1}{2}(-K_{V})^{3} + 1$) is
called (anticanonical) \emph{degree} (respectively, \emph{genus})
of $V$.

\smallskip

\item For a Weil divisor $D$ on a normal algebraic variety
$V$, we denote by $\mathcal{O}_{V}(D)$ the corresponding
divisorial sheaf on $V$ (sometimes we denote both by
$\mathcal{O}_{V}(D)$ (or by $D$)).

\smallskip

\item For a vector bundle $\mathcal{E}$ on smooth projective variety
$V$, we denote by $c_{i}(\mathcal{E})$ the $i$-th Chern class of
$\mathcal{E}$.

\smallskip

\item We denote by $T_{P}(V)$ the Zariski tangent space to an algebraic variety $V$ at a point $P \in
V$. For $V$ smooth and a smooth hypersurface $D \subset V$, we
denote by $T_{V}\langle D \rangle$ the subsheaf of the tangent
sheaf on $V$ which consists of all vector fields tangent to $D$.

\smallskip

\item For a Cartier divisor $M$ on a normal projective variety $V$, we denote by $|M|$ the corresponding complete linear system on $V$.
For an algebraic cycle $Z$ on $V$, we denote by $|M - Z|$ the
linear subsystem in $|M|$ which consists of all divisors passing
through $Z$. For a linear system $\mathcal{M}$ on $V$ without base
components, we denote by $\Phi_{\mathcal{M}}$ the corresponding
rational map.

\smallskip

\item For a birational map $\psi: V' \dashrightarrow V$ between normal projective varieties
and an algebraic cycle $Z$ (respectively, a linear system
$\mathcal{M}$) on $V$, we denote by $\psi_{*}^{-1}(Z)$
(respectively, by $\psi_{*}^{-1}(\mathcal{M})$) the proper
transform of $Z$ (respectively, of $\mathcal{M}$) on $V'$.

\smallskip

\item We denote by $\mathbb{F}_{n}$ the Hirzebruch
surface with the class of a fiber $l$ and the minimal section $h$
of the natural projection $\mathbb{F}_{n} \to \mathbb{P}^1$ such
that $(h^{2}) = -n$, $n \in \mathbb{Z}_{\geqslant 0}$.

\end{itemize}

\section{Preliminaries}
\label{section:preliminaries}

In what follows, $X$ is a Fano $3$-fold of genus $36$ (or degree
$70$). Let us present the construction and some properties of $X$
(see \cite{karz-1} for more details).

Consider the weighted projective space $\mathbb{P} :=
\mathbb{P}(1,1,4,6)$ with weighted homogeneous coordinates
$x_{0}$, $x_{1}$, $x_{2}$, $x_{3}$ of weights $1$, $1$, $4$, $6$,
respectively. $\mathbb{P}$ is a Fano $3$-fold of degree $72$.
Furthermore, the linear system $|-K_{\mathbb{P}}|$ gives an
embedding of $\mathbb{P}$ in $\mathbb{P}^{38}$ such that the image
$\Phi_{\scriptscriptstyle|-K_{\mathbb{P}}|}(\mathbb{P})$ is an
intersection of quadrics. In what follows, we assume that
$\mathbb{P} \subset \mathbb{P}^{38}$ is anticanonically embedded.
Then $L := \mathrm{Sing}(\mathbb{P})$ is a line on $\mathbb{P}$
with respect to this embedding. Moreover, there are two points $P$
and $Q$ on $L$ such that the singularities $P \in \mathbb{P}$, $Q
\in\mathbb{P}$ are of types $\displaystyle\frac{1}{6}(4,1,1)$,
$\displaystyle\frac{1}{4}(2,1,1)$, respectively, and for every
point $O \in L \setminus\{P,~Q \}$ the singularity $O \in
\mathbb{P}$ is analytically isomorphic to $(0, o) \in \mathbb{C}
\times W$, where $o \in W$ is the singularity of type
$\displaystyle\frac{1}{2}(1,1)$ (see \cite[Example 2.13]{karz-1}).

\begin{proposition}
\label{theorem:unique-line} $L$ is the unique line on
$\mathbb{P}$.
\end{proposition}

\begin{proof}
Let $L_{0} \ne L$ be another line on $\mathbb{P}$. Since
$-K_{\mathbb{P}} \sim \mathcal{O}_{\mathbb{P}}(12)$, we have
\begin{equation}
\label{bad-intersection} \mathcal{O}_{\mathbb{P}}(1) \cdot L_0 =
\frac{1}{12},
\end{equation}
which implies that $L \cap L_0 \ne \emptyset$. Consider the
crepant resolution $\phi: T \longrightarrow \mathbb{P}$ of
$\mathbb{P}$. Set $L'_0 := \phi_{*}^{-1}(L_0)$, $E_Q :=
\phi^{-1}(Q)$, $E_P := \phi^{-1}(P)$ and $E_L :=
\overline{\phi^{-1}(L \setminus{\{P,Q\}})}$, the Zariski closure
in $T$ of $\phi^{-1}(L \setminus{\{P,Q\}})$. These are all the
components of the $\phi$-exceptional locus. Furthermore, we have
$E_P = E^{(1)}_P \cup E^{(2)}_P$, where $E^{(i)}_P$ are
irreducible components of the divisor $E_P$ such that $E^{(1)}_P
\cap E_L = \emptyset$ and $E^{(2)}_P \cap E_L \ne \emptyset$ (see
\cite[Example 2.13]{karz-1} for the explicit construction of
$\phi$).

Since $\rho(\mathbb{P}) = 1$, the group $N_{1}(T)$ is generated by
the classes of $\phi$-exceptional curves and some curve $Z$ on $T$
such that $R: = \mathbb{R}_{+}[Z]$ is the $K_{T}$-negative
extremal ray (see \cite[Lemmas 4.2, 4.3]{prok}). In particular,
since $-K_{T} \cdot L'_0 = 1$, \cite[Lemmas 4.2, 4.3]{prok}
implies that
\begin{equation}
\label{numerical-equality-for-L-0} L'_0 \equiv Z + E^*,
\end{equation}
where $E^*$ is a linear combination with nonnegative coefficients
of irreducible $\phi$-exceptional curves. Further, the linear
projection $\pi_L$ of $\mathbb{P}$ from $L$ is given by the linear
system $\mathcal{H} \subset |-K_{\mathbb{P}}|$ of all hyperplane
sections of $\mathbb{P}$ containing $L$. In addition, $\pi_L$ maps
$L_0$ to the point because $L \cap L_0 \ne \emptyset$ and
$\mathbb{P}$ is the intersection of quadrics. On the other hand,
$\phi$ factors through the blow up of $\mathbb{P}$ at $L$ (see
\cite{karz-1}, \cite{karz-2}). Hence the linear system
$\phi_{*}^{-1}\mathcal{H}$ is basepoint-free on $T$ and $H \cdot
L'_0 = 0$ for $H \in \phi_{*}^{-1}\mathcal{H}$. In particular, $H
\in |-K_T - E_L|$.

\begin{lemma}
\label{theorem:few-components} In
\eqref{numerical-equality-for-L-0}, the support
$\mathrm{Supp}(E^*)$ of $E^*$ is either $\emptyset$ or $e_P$,
where $e_P \subset E^{(1)}_P$.
\end{lemma}

\begin{proof}
As we saw, the face of the Mori cone $\overline{NE}(T)$, which
corresponds to the nef divisor $H$, contains the class of the
curve $L'_0$. Then from \eqref{numerical-equality-for-L-0} we get
$$
H \cdot Z = H \cdot E^* = 0.
$$
In particular, $H$ intersects trivially every curve in
$\mathrm{Supp}(E^*)$. On the other hand, we have
$\mathrm{Supp}(E^*) \subseteq \{e_{P}, e_{Q}, e_{L}\}$, where
$e_P$, $e_Q$, $e_L$ are the curves in $E_P$, $E_Q$, $E_L$,
respectively. But for $e_P \subset E^{(2)}_P$ intersections $H
\cdot e_P$, $H \cdot e_Q$, $H \cdot e_L$ are all non-zero. Thus,
$\mathrm{Supp}(E^*)$ is either $\emptyset$ or $e_P$, where $e_P
\subset E^{(1)}_P$.
\end{proof}

Consider the extremal contraction $f_{R} : T \longrightarrow T'$
of $R$. The morphism $f_{R}$ is birational with the exceptional
divisor $E_R$ (see \cite{karz-1}, \cite{karz-2}).

\begin{lemma}
\label{theorem:non-nef} The divisor $-K_{T'}$ is not nef.
\end{lemma}

\begin{proof}
Suppose that $-K_{T'}$ is nef, i.e., $T'$ is a weak Fano $3$-fold
(with possibly non-Gorenstein singularities). If $T'$ has only
terminal factorial singularities, then since $(-K_{T'})^{3}
\geqslant (-K_{T})^{3} = 72$ (see \cite[Proposition-definition
4.5]{prok}), $T'$ is a terminal $\mathbb{Q}$-factorial
modification either of $\mathbb{P}(1,1,1,3)$ or of
$\mathbb{P}(1,1,4,6)$. In particular, either $\rho(Y') = 5$ or
$\rho(Y') = 2$ (see \cite{karz-1}, \cite{karz-2}). On the other
hand, $\rho(T') = \rho(T) - 1 = 4$, a contradiction.

Thus, the singularities of $T'$ are worse than factorial. In this
case, $f_{R}(E_{R})$ is a point (see \cite[Proposition-definition
4.5]{prok}) and we get
\begin{equation}
\label{non-empty-E-R-2} E_P \cap E_R = E_Q \cap E_R = \emptyset.
\end{equation}
On the other hand, it follows from
\eqref{numerical-equality-for-L-0} that $-K_{\mathbb{P}} \cdot
\phi_{*}(Z) = 1$, i.e., $\phi(Z)$ is a line on $\mathbb{P}$. In
particular, as for $L_0$ above, we have $\phi_{*}(Z) \cap L \ne
\emptyset$. But then \eqref{non-empty-E-R-2} implies that $0 =
K_{T} \cdot Z = -1$, a contradiction.
\end{proof}

It follows from Lemma~\ref{theorem:non-nef} that $E_R =
\mathbb{F}_{1}$ or $\mathbb{P}^{1} \times \mathbb{P}^{1}$ (see
\cite[Proposition-definition 4.5]{prok}). But if $E_R  =
\mathbb{F}_{1}$, then $\phi(E_R)$ is a plane on $\mathbb{P}$ such
that $L \not \subset \phi(E_R)$ (see \cite[Proposition-definition
4.5]{prok}). This implies that there is a line on $\mathbb{P}$ not
intersecting $L$, a contradiction (see \eqref{bad-intersection}).
Finally, in the case when $E_R = \mathbb{P}^{1} \times
\mathbb{P}^{1}$, we have $Z \subset E_R = E_L$ (see
\cite[Proposition-definition 4.5]{prok}), and if
$\mathrm{Supp}(E^*) = \emptyset$ in
\eqref{numerical-equality-for-L-0}, then $L_0 = L$, a
contradiction. Hence, by Lemma~\ref{theorem:few-components}, we
get $\mathrm{Supp}(E^*) = e_P$, where $e_P \subset E^{(1)}_P$.
Further, on $E_R$ we have:
$$
Z \sim l, \qquad E_P\big\vert_{E_R} = E^{(2)}_P\big\vert_{E_R}
\sim h \sim E_Q\big\vert_{E_R},
$$
which implies that $E^{(2)}_P \cdot Z = E_Q \cdot Z = 1$. On the
other hand, since $L_0 \ne L$, we have either $E^{(2)}_P \cdot
L'_0 = 0$ or $E_Q \cdot L'_0 = 0$. Then, intersecting
\eqref{numerical-equality-for-L-0} with $E^{(2)}_P$ and $E_Q$, we
get a contradiction because $E^{(2)}_P \cdot e_P$ and $E_Q \cdot
e_P \geqslant 0$.

Thus, we get $L_0 = L$, a contradiction.
Proposition~\ref{theorem:unique-line} is completely proved.
\end{proof}

Coming back to the construction of $X$, take any point $O$ in $L
\setminus\{P,~Q \}$ and consider the linear projection $\pi:
\mathbb{P} \dashrightarrow \mathbb{P}^{37}$ from $O$. Then the
image of $\pi$ is a Fano $3$-fold $X_O$ of degree $70$ (see
\cite{karz-1}, \cite{karz-2}).

\begin{proposition}
\label{theorem:unique-up-to-isomorphism} For any point $O'$ in $L
\setminus\{P,~Q,~O\}$, the image of the linear projection
$\mathbb{P} \dashrightarrow \mathbb{P}^{37}$ from $O'$ is a Fano
$3$-fold $X_{O'}$ isomorphic to $X_O$.
\end{proposition}

\begin{proof}
In the above notation, $L$ is given by equations $x_0 = x_1 = 0$
on $\mathbb{P}$, with equations of $P$ and $Q$ being $x_0 = x_1 =
x_2 = 0$ and $x_0 = x_1 = x_3 = 0$, respectively (see
\cite[5.15]{iano-fletcher}). Then the torus $(\mathbb{C}^{*})^3$,
acting on $\mathbb{P}$, acts transitively on the set $L
\setminus\{P,~Q \}$, which induces an isomorphism $X_{O'} \simeq
X_O$.
\end{proof}

In what follows, because of
Proposition~\ref{theorem:unique-up-to-isomorphism}, we fix the
point $O \in L \setminus\{P,~Q \}$, the linear projection $\pi:
\mathbb{P} \dashrightarrow \mathbb{P}^{37}$ from $O$, and denote
the image of $\pi$ by $X$. Let us construct a terminal
$\mathbb{Q}$-factorial modification of $X$. Consider the blow up
$\sigma: W \longrightarrow \mathbb{P}$ of $\mathbb{P}$ at $O$, and
the following commutative diagram:
$$
\xymatrix{
&&W\ar@{->}[ld]_{\sigma}\ar@{->}[rd]^{\mu}&&\\%
&\mathbb{P}\ar@{-->}[rr]_{\pi}&&X.&}
$$
The type of the singularity $O \in \mathbb{P}$ implies that $W$
has at most canonical Gorenstein singularities. Moreover, we have
$\mathrm{Sing}(W) = \sigma_{*}^{-1}(L)$ and the singularities of
$W$ are exactly of the same kind as of $\mathbb{P}$, i.e., locally
near every point in $\mathrm{Sing}(W)$, $W$ is isomorphic to
$\mathbb{P}$. Then, resolving the singularities of $W$ in the same
way as for $\mathbb{P}$, we arrive at the birational morphism
$\tau : Y \longrightarrow W$, with $Y$ being smooth and $K_Y =
\tau^{*}(K_W)$ (see \cite{karz-1}, \cite{karz-2}). Set $f := \tau
\circ \mu : Y \longrightarrow X$.

\begin{proposition}
\label{theorem:uniqueness-up-to-isomorphism} $f : Y
\longrightarrow X$ is a terminal $\mathbb{Q}$-factorial
modification of $X$. Moreover, $Y$ is unique up to isomorphism,
i.e., every smooth weak Fano $3$-fold of degree $70$ is isomorphic
to $Y$.
\end{proposition}

\begin{proof}
The linear projection $\pi$ is given by the linear system
$\mathcal{H} \subset |-K_{\mathbb{P}}|$ of all hyperplane sections
of $\mathbb{P}$ passing through $O$. For a general $H \in
\mathcal{H}$, we have
$$
\sigma_{*}^{-1}(H) = \sigma^{*}(H) - E_{\sigma},
$$
where $E_{\sigma}$ is the $\sigma$-exceptional divisor. On the
other hand, from the adjunction formula we get
$$
K_{W} = \sigma^{*}(K_{\mathbb{P}}) + E_{\sigma}.
$$
Thus, the morphism $\mu: W \longrightarrow X$ is given by the
linear system $\sigma_{*}^{-1}(\mathcal{H}) \subseteq |-K_{W}|$.
Furthermore, since $\mathbb{P}$ is an intersection of quadrics,
$\pi$ is a birational map, which implies that $\mu$ and $f$ are
also birational with $K_Y = f^{*}(K_{X})$. In particular,
$(-K_{Y})^{3} = (-K_{X})^{3} = 70$.

Thus, it remains to prove that every smooth weak Fano $3$-fold of
degree $70$ is isomorphic to $Y$. Let $Y'$ be another smooth weak
Fano $3$-fold of degree $70$. Then its image under the morphism
$f' := \Phi_{\scriptscriptstyle|-nK_{Y'}|}$, $n \in \mathbb{N}$,
is a Fano threefold $X'$ such that $K_{Y'} \equiv f'^{*}(K_{X'})$
(see \cite{kawamata}). Since $(-K_{Y'})^3 = (-K_{X'})^3 = 70$, we
get $X' \simeq X$ and $Y'$ is a terminal $\mathbb{Q}$-factorial
modification of $X$. Then, since $Y$ and $Y'$ are relative minimal
models over $X$, the induced birational map $Y \dashrightarrow Y'$
is either an isomorphism or a sequence of $K_{Y}$-flops over $X$
(see \cite{kollar-mori}).

\begin{lemma}
\label{theorem:y-is-unique} Every $K_{Y}$-trivial extremal
birational contraction $f_1 : Y \longrightarrow Y_1$ is
divisorial.
\end{lemma}

\begin{proof}
Suppose that $f_1$ is small. In the notation from the proof of
Proposition~\ref{theorem:unique-line}, denote by $E_{Y,L}$,
$E_{P,L}^{(i)}$, $E_{Q,L}$ the proper transforms on $Y$ of
surfaces $E_L$, $E_{P}^{(i)}$, $E_{Q}$, respectively. The
resolution $\tau : Y \longrightarrow W$ (or $\phi: T
\longrightarrow \mathbb{P}$) is locally toric near
$\mathrm{Sing}(W)$. In particular, we have $E_{P,L}^{(1)} \simeq
\mathbb{F}_4$, $E_{P,L}^{(2)} \simeq \mathbb{F}_2$, $E_{Q,L}
\simeq \mathbb{F}_2$, $E_{Y,L} \simeq \mathbb{F}_m$ for some $m
\in \mathbb{N}$ (see \cite[Example 2.13]{karz-1}), and hence the
only possibility for $f_1$ is to contract the curve $Z = h$ on
$E_{Y,L}$ such that $\tau(Z) = \sigma_{*}^{-1}(L)$.

On the other hand, $Y$ is obtained by the blow up of the $3$-fold
$T$ at the curve $\phi^{-1}(O) \simeq \mathbb{P}^1$ (see
\cite{karz-1}, \cite{karz-2}). Furthermore, since $\mathbb{P}$ is
singular along the line, we have $E_L \simeq \mathbb{P}^1 \times
\mathbb{P}^1$ (see \cite[Proposition-definition 4.5]{prok}), and
hence $E_{Y,L} \simeq \mathbb{P}^1 \times \mathbb{P}^1$, a
contradiction.
\end{proof}

It follows from Lemma~\ref{theorem:y-is-unique} that $Y' \simeq
Y$. Proposition~\ref{theorem:uniqueness-up-to-isomorphism} is
completely proved.
\end{proof}

\begin{corollary}
\label{theorem:singularities-of-x} $\mathrm{Sing}(X)$ consists of
a unique point.
\end{corollary}

\begin{proof}
Since the morphism $\mu: W \longrightarrow X$ is given by the
linear system $\sigma_{*}^{-1}(\mathcal{H}) \subseteq |-K_{W}| =
|\sigma^{*}(-K_{\mathbb{P}}) - E_{\sigma}|$ (see the proof of
Proposition~\ref{theorem:uniqueness-up-to-isomorphism}), it
contracts only $\sigma_{*}^{-1}(L) = \mathrm{Sing}(W)$ to the
unique singular point on $X$ (see
Proposition~\ref{theorem:unique-line}).
\end{proof}

\begin{corollary}
\label{theorem:k-x-generates-pic-of-x} We have $\mathrm{Pic}(X) =
\mathbb{Z} \cdot K_X$ and $\mathrm{Cl}(X) = \mathbb{Z} \cdot K_X
\oplus \mathbb{Z} \cdot E$, where $E := \mu_{*}(E_{\sigma})$.
\end{corollary}

\begin{proof}
This follows from the construction of $X$ and equalities
$\rho(\mathbb{P}) = 1$, $(-K_{X})^3 = 70$.
\end{proof}

\begin{remark}
\label{remark:x-is-anti-canonically-embedded} It follows from the
construction of $X$ that $f = \Phi_{\scriptscriptstyle|-K_{Y}|}$
and $X \subseteq \mathbb{P}^{37}$ is anticanonically embedded.
\end{remark}

\begin{remark}
\label{remark:picard-and-algebraic-cycles} Since $Y$ is a smooth
weak Fano $3$-fold, we have $\mathrm{Pic}(Y) \simeq H^{2}(Y,
\mathbb{Z})$ (see \cite[Proposition 2.1.2]{isk-prok}) and
$H^{2}(Y, \mathcal{O}_{Y}) = 0$ by Kawamata--Viehweg Vanishing
Theorem.
\end{remark}

It follows from Corollary~\ref{theorem:singularities-of-x} that a
general surface $S \in |-K_{X}|$ is smooth. Furthermore,
Corollary~\ref{theorem:k-x-generates-pic-of-x} implies that the
cycles $K_{X} \big\vert_{S}$ and $E \big\vert_{S}$ are not
divisible in $\mathrm{Pic}(S)$, linearly independent in $H^{2}(S,
\mathbb{Q})$, and hence they generate a primite sublattice $R_S$
in $\mathrm{Pic}(S)$. It follows from the construction of $X$ that
all lattices $R_S$, $S \in |-K_{X}|$, are isomorphic to the
lattice $\mathrm{R} \simeq \mathbb{Z}^2$ with the associated
quadratic form $70x^2 + 4xy - 2y^2$, and we can consider the
moduli stack $\mathcal{K} := \mathcal{K}_{36}^{\mathrm{R}}$ of
$\mathrm{K3}$ surfaces of type $\mathrm{R}$ (see
\cite[(2.3)]{beauville}). $\mathcal{K}$ is actually an algebraic
space because the forgetful map $\mathcal{K} \longrightarrow
\mathcal{K}_{36}$ is representable and $1$-to-$1$ in our case (see
\cite[(2.5)]{beauville}).\footnote{It can be also easily seen that
the class of a $(-2)$-curve in $\mathrm{Pic}(S)$ is unique and
generated by the conic $E\vert_{S}$.}

\begin{proposition}[see \cite{beauville}]
\label{theorem:moduli-space-of-polarized-k-3-surfaces} Let $S$ be
the $\mathrm{K3}$ surface of type $\mathrm{R}$. Then

\begin{itemize}

\item the first order deformations of $(S,
\mathrm{R})$ are parameterized by the orthogonal of
$c_{1}(\mathrm{R}) \subset H^{1}(S, \Omega_{S}^{1})$ in $H^{1}(S,
T_{S})$;

\item the space $\mathcal{K}$ is smooth, irreducible, of
dimension $18$.

\end{itemize}

\end{proposition}

\section{Proof of Theorem~\ref{theorem:main}}
\label{section:proof}

We use the notation and conventions of
Section~\ref{section:preliminaries}. Since $f : Y \longrightarrow
X$ is the crepant resolution (see
Proposition~\ref{theorem:uniqueness-up-to-isomorphism}), it
follows from Corollary~\ref{theorem:singularities-of-x} that we
can assume a general $S \in |-K_{X}|$ to be a surface in
$|-K_{Y}|$ on $Y$. We can also assume that $S \cap \mathrm{Exc}(f)
= \emptyset$ for the $f$-exceptional locus $\mathrm{Exc}(f)$.
Further, it follows from
Remark~\ref{remark:x-is-anti-canonically-embedded} that the points
in $(\mathbb{P}^{37})^{*}$, corresponding to such $S$'s, form an
open subset $U \subset (\mathbb{P}^{37})^{*}$. Consider the
natural (faithful) action of the group $G := \mathrm{Aut}(Y)$ on
$U$. Shrinking $U$ if necessary, we obtain the following

\begin{proposition}
\label{theorem:locally-free-action} The geometric quotient $U/G$
exists as a smooth scheme.
\end{proposition}

\begin{proof}
Let us calculate the group $G$ first. Take $g \in
\mathrm{Aut}(\mathbb{P})$ to be an automorphism of $\mathbb{P}$
which fixes the point $O$. Then $g$ lifts to the automorphism of
$Y$ (see the construction of $X$ and $Y$ in
Section~\ref{section:preliminaries}). Conversely, take any $g \in
G$.

\begin{lemma}
\label{theorem:tau-is-g-equivariant} The morphism $\tau: Y
\longrightarrow W$ is $g$-equivariant.
\end{lemma}

\begin{proof}
Since the morphism $f = \Phi_{\scriptscriptstyle|-K_{Y}|} : Y
\longrightarrow X$ is $g$-equivariant (see
Remark~\ref{remark:x-is-anti-canonically-embedded}), it follows
from the construction of $Y$ in
Section~\ref{section:preliminaries} that the irreducible
components of $\mathrm{Exc}(f)$ are all $g$-invariant. Thus, since
$\mathrm{Pic}(Y)$ is generated by $K_Y$, the irreducible
components of $E_f$ and $E_{Y,\sigma} :=
\tau_{*}^{-1}(E_{\sigma})$, it is enough to prove that
$g(E_{Y,\sigma}) = E_{Y,\sigma}$. Suppose that $g(E_{Y,\sigma})
\ne E_{Y,\sigma}$. Then, since all the curves in $E_{\sigma}$
(respectively, in $\tau_{*}(g(E_{Y,\sigma}))$) are numerically
proportional and $\tau$ is divisorial, we must have $E_{\sigma}
\cap \tau_{*}(g(E_{Y,\sigma})) = \emptyset$. The latter implies
that there exists a curve $C \equiv \sigma_{*}(-K_{W} \cdot
\tau_{*}(g(E_{Y,\sigma})))$ on $\mathbb{P}$ with $-K_{\mathbb{P}}
\cdot C = 4$ and $C \cap L = \emptyset$. On the other hand, since
$-K_{\mathbb{P}} \sim \mathcal{O}_{\mathbb{P}}(12)$, we get
$\mathcal{O}_{\mathbb{P}}(1) \cdot C = \displaystyle\frac{1}{3}$,
a contradiction.
\end{proof}

It follows from Lemma~\ref{theorem:tau-is-g-equivariant} that $g$
acts on $W$. Further, considering the induced $g$-action on the
cone $\overline{NE}(W)$, we obtain, since $\mathrm{Pic}(W) =
\mathbb{Z} \cdot K_W \oplus \mathbb{Z} \cdot E_{\sigma}$, that
$\sigma: W \longrightarrow \mathbb{P}$ is $g$-equivariant. The
latter gives a $g$-action on $\mathbb{P}$ with the fixed point
$O$.

Thus, $G$ is isomorphic to the stabilizer in
$\mathrm{Aut}(\mathbb{P})$ of the point $O$, and to describe the
$G$-action on $U$ we may consider the action of the corresponding
subgroup in $\mathrm{Aut}(\mathbb{P})$ on the linear system
$|-K_{\mathbb{P}} - O|$. Note that, since $P \in \mathbb{P}$, $Q
\in \mathbb{P}$, $O \in \mathbb{P}$ are the pairwise
non-isomorphic singularities, every $g \in G$ fixes every point on
$L$. Finally, since $\mathcal{O}_{\mathbb{P}}(1)$,
$\mathcal{O}_{\mathbb{P}}(4)$, $\mathcal{O}_{\mathbb{P}}(6)$ are
$G$-invariant, the $g$-action on $\mathbb{P}$ can be described as
follows:
\begin{eqnarray}
\label{g-action} x_0 \mapsto ax_0 + bx_1, \\
\nonumber x_1 \mapsto cx_0 + dx_1, \\
\nonumber x_2 \mapsto \lambda^{4} x_2 + f_{4}(x_{0},x_{1}),\\
\nonumber x_3 \mapsto \lambda^{6} x_3 + x_{2}f_{2}(x_{0},x_{1}) +
f_{6}(x_{0},x_{1}),
\end{eqnarray}
where $\lambda \in \mathbb{C}^*$, $\bigl(\begin{smallmatrix} a & b \\
c & d\end{smallmatrix}\bigr) \in GL(2, \mathbb{C})/\{\pm 1\}$,
$f_{i} := f_{i}(x_{0}, x_{1})$ are arbitrary homogeneous
polynomials of degree $i$ in $x_{0}$, $x_{1}$. On the other hand,
since $-K_{\mathbb{P}} \sim \mathcal{O}_{\mathbb{P}}(12)$, a
general element in $|-K_{\mathbb{P}} - O|$ can be given by the
equation
\begin{eqnarray}
\label{s-equation} \alpha x_{3}^{2} + x_{2}^3 + a_{6}(x_{0},
x_{1})x_{3} + a_{2}(x_{0}, x_{1})x_{2}x_{3} + a_{4}(x_{0},
x_{1})x_{2}^{2} + a_{8}(x_{0}, x_{1})x_{2} + a_{12}(x_{0}, x_{1})
= 0
\end{eqnarray}
on $\mathbb{P}$, where $a_{i} := a_{i}(x_{0}, x_{1})$ are
arbitrary general homogeneous polynomials in $x_{0}$, $x_{1}$ of
degree $i$, and $\alpha \in \mathbb{C}^*$ is fixed.

Take a general surface $S_0$ on $\mathbb{P}$ with the equation
\eqref{s-equation} such that $a_2 = a_4 = a_6 = 0$.

\begin{lemma}
\label{theorem:generic-surface-not-g-invariant} If $S_0$ is
$g$-invariant for some $g \ne \mathrm{id}$ from \eqref{g-action},
then $f_2 = f_4 = f_6 = 0$, $c = b = 0$, $a = d = \sqrt{-1}$,
$\lambda^4 = 1$.
\end{lemma}

\begin{proof}
$g$-invariance of $S_0$ implies that $f_2 = f_4 = f_6 = 0$ and
\begin{eqnarray}
\label{not-free-g-action} a_{8}(x_{0}, x_{1}) = a_{8}(ax_0 + bx_1,
cx_0 + dx_1), \qquad a_{12}(x_{0}, x_{1}) = a_{12}(ax_0 + bx_1,
cx_0 + dx_1).
\end{eqnarray}
Without loss of generality we may assume that $a_8 =
x_{0}x_{1}b_6$ for some $b_{6} := b_{6}(x_{0}, x_{1})$ coprime
to $x_0$ and $x_1$. Then \eqref{not-free-g-action} and generality of $S_0$ imply that $\bigl(\begin{smallmatrix} a & b \\
c & d\end{smallmatrix}\bigr) = \bigl(\begin{smallmatrix} a & 0 \\
0 & d \end{smallmatrix}\bigr)$, and we get:
$$
a^{12} = 1, \qquad a^{i+1}d^{7-i} = 1, \qquad a^{i}d^{6-i} =
a^{j}d^{6-j}
$$
for all $0 \leqslant i, \ j \leqslant 6$. In particular, $a = d$,
$a^8 = a^{12} = 1$, i.e., $a = d = \sqrt{-1}$. Finally, since $x_2
\mapsto \lambda^4 x_2$ (see \eqref{g-action}) and hence
$a_{8}(x_{0}, x_{1}) = \lambda^4 a_{8}(x_{0}, x_{1})$ (see
\eqref{s-equation}), we get $\lambda^4 = 1$.
\end{proof}

\begin{lemma}
\label{theorem:when-we-get-identity} Let $g \in G$, given by
\eqref{g-action}, be such that $f_2 = f_4 = f_6 = 0$, $c = b = 0$,
$a = d = \sqrt{-1}$, $\lambda = \pm \sqrt{-1}$. Then $g =
\mathrm{id}$.
\end{lemma}

\begin{proof}
We have
$$
g([x_{0}:x_{1}:x_{2}:x_{3}]) = [\sqrt{-1}
x_{0}:\sqrt{-1}x_{1}:(\sqrt{-1})^{4}x_{2}:(\sqrt{-1})^{6}x_{3}] =
[x_{0}:x_{1}:x_{2}:x_{3}]
$$
on $\mathbb{P}$. Hence $g = \mathrm{id}$.
\end{proof}

It follows from
Lemmas~\ref{theorem:generic-surface-not-g-invariant} and
\ref{theorem:when-we-get-identity}, since $\lambda^4 = 1$ implies
$\lambda^2 = \pm 1$, that the stabilizer of $S_0$ in $G$ is a
group of order $2$, generated by some $g_0 \in G$ with $\lambda^2
= 1$ (see \eqref{g-action}). Consider the normal algebraic
subgroup $G' \subset G$ generated by $g^{-1}g_{0}g$ for all $g \in
G$, i.e., generators of $G'$ are all the elements in $G$ for which
$f_4 = 0$, $c = b = 0$, $a = d = \sqrt{-1}$ and $\lambda = 1$ in
\eqref{g-action}. Then the $G'$-action on $U$ is proper, and we
can consider the geometric quotient $U' := U/G'$, which exists as
a normal scheme (see \cite{popp}). Further, take the $G'' :=
G/G'$-equivariant factorization map $\pi_{G} : U \longrightarrow
U'$ and consider the induced $G''$-action on $U'$. Shrinking $U$
if necessary, we obtain

\begin{lemma}
\label{theorem:free-action-near-l-0} The $G''$-action on $U'$ is
free.
\end{lemma}

\begin{proof}
Let $S'_0$ be the image on $U'$ of $S_0$ under $\pi_{G}$. Then we
have $G'' \cdot S'_0 \simeq G''$ for the $G''$-orbit of $S'_0$,
and, by the dimension count, there exists a Zariski open subset in
$U'$ with a free $G''$-action.
\end{proof}

Lemma~\ref{theorem:free-action-near-l-0} and \cite{popp} imply
that the geometric quotient $U/G \simeq U'/G''$ exists as a smooth
scheme. Proposition~\ref{theorem:locally-free-action} is
completely proved.
\end{proof}

Set $\mathcal{F} := U/G$ to be the scheme from
Proposition~\ref{theorem:locally-free-action}. It follows from
Proposition~\ref{theorem:uniqueness-up-to-isomorphism} and
Remark~\ref{remark:x-is-anti-canonically-embedded} that
$\mathcal{F}$ is a (coarse) moduli space which parameterizes the
pairs $(Y^{\sharp}, S^{\sharp})$ consisting of smooth weak Fano
$3$-fold $Y^{\sharp}$ of degree $70$ and smooth surface
$S^{\sharp} \in |-K_{Y^{\sharp}}|$ (cf. \cite[(2.2)]{beauville}).
These give the following

\begin{lemma}
\label{theorem:tangent-space-is-H-1} For $o := (Y, S) \in
\mathcal{F}$, we have $H^{1}(Y, T_{Y}\langle S \rangle) =
T_{o}\mathcal{F}$.
\end{lemma}

\begin{proof}
This follows from the fact that $\mathcal{F}$ is smooth and
$H^{1}(Y, T_{Y}\langle S \rangle)$ parameterizes the first order
deformations of $(Y, S)$ (see \cite[Proposition 1.1]{beauville}).
\end{proof}

Consider the forgetful morphism $s: \mathcal{F} \longrightarrow
\mathcal{K}$, which sends $(Y, S)$ to $S$.

\begin{lemma}
\label{theorem:s-is-surjective} $s$ is generically surjective.
\end{lemma}

\begin{proof}
Consider the restriction map $r: T_{Y}\langle S \rangle
\longrightarrow T_S$. It fits into the exact sequence
\begin{equation}
\label{eq-1} 0 \longrightarrow \Omega_{Y}^{2} \longrightarrow
T_{Y}\langle S \rangle \stackrel{r}{\longrightarrow} T_S
\longrightarrow 0,
\end{equation}
since $\mathrm{Ker}(r) = T_{Y}(-S)$ is a subsheaf of $T_{Y}\langle
S \rangle$ consisting of the vector fields vanishing along $S$,
for which we have $T_{Y}(-S) \simeq \Omega_{Y}^{2}$. From
\eqref{eq-1} we get the exact sequence
$$
H^{1}(Y, T_{Y}\langle S \rangle)
\stackrel{H^{1}(r)}{\longrightarrow} H^{1}(S, T_{S})
\stackrel{\partial}{\longrightarrow} H^{2}(Y, \Omega_{Y}^{2}).
$$
The map $\partial$ is dual to the restriction map $i: H^{1}(Y,
\Omega_{Y}^{1}) \longrightarrow H^{1}(S, \Omega_{S}^{1})$ (see
\cite{beauville}). In particular, $\mathrm{Ker}(\partial)$ is the
orthogonal of $\mathrm{Im}(i)$. On the other hand, we have
$\mathrm{Im}(i) = \mathbb{Z} \cdot c_{1}(K_{Y} \big\vert_{S})
\oplus \mathbb{Z} \cdot c_{1}(\tau_{*}^{-1}(E_{\sigma})
\big\vert_{S}) \simeq \mathbb{Z} \cdot K_{X} \big\vert_{S} \oplus
\mathbb{Z} \cdot E \big\vert_{S}$ (see
Corollary~\ref{theorem:k-x-generates-pic-of-x} and
Remark~\ref{remark:picard-and-algebraic-cycles}), and hence
$H^{1}(r)$ coincides with the tangent map to $s$ at $(Y, S)$, with
$\mathrm{Im}(H^{1}(r)) = \mathrm{Ker}(\partial)$ being the tangent
space to $\mathcal{K}$ at $S$ (see
Lemma~\ref{theorem:tangent-space-is-H-1} and
Proposition~\ref{theorem:moduli-space-of-polarized-k-3-surfaces}).
Thus, since $\mathcal{K}$ is irreducible (see
Proposition~\ref{theorem:moduli-space-of-polarized-k-3-surfaces}),
we get that $s$ is generically surjective.
\end{proof}

Theorem~\ref{theorem:main} is completely proved.

\end{document}